\input amstex
\documentstyle{amsppt}
\magnification=\magstep1 \TagsOnRight
\NoBlackBoxes

\hoffset1 true pc
\voffset2 true pc
\hsize36 true pc
\vsize52 true pc

\tolerance=2000
\def\m1{^{-1}}
\def\ov1{\overline}
\def\gp#1{\langle#1\rangle}

\def\NP{\left(\smallmatrix
 E_{n}&0\\
\endsmallmatrix\right)}

\def\NL{\left(\smallmatrix
 0&E_{n}\\
\endsmallmatrix\right)}

\catcode`\@=11
\def\logo@{}
\catcode`\@=\active

\topmatter
\title
Indecomposable  linear groups
\endtitle
\author
V.A.~Bovdi, V.P.~Rudko
\endauthor
\dedicatory Dedicated to the memory of Professor A.V.~Roiter
\enddedicatory
\leftheadtext\nofrills{ V.A.~Bovdi, V.P.~Rud'ko } \rightheadtext
\nofrills { Indecomposable  linear groups } \abstract Let $G$ be a
noncyclic group of order $4$, and let  $K=\Bbb Z$, $\Bbb Z_{(2)}$
and $\Bbb Z_{2}$ be the ring of rational integers, the
localization of $\Bbb Z$ at the prime $2$ and the ring of $2$-adic
integers, respectively. We describe, up to conjugacy, all of the
indecomposable subgroups in the group $GL(m,K)$  which are
isomorphic to $G$.
\endabstract
\thanks
The research was supported by OTKA  No.T 037202 and  No.T 038059
\endthanks

\address
\hskip-\parindent  Institute of Mathematics, University of
Debrecen\newline P.O.  Box 12, H-4010 Debrecen, Hungary\newline
Institute of Mathematics and Informatics, College of Ny\'\i
regyh\'aza\newline S\'ost\'oi \'ut 31/b, H-4410 Ny\'\i regyh\'aza,
Hungary
\newline
E-mail: vbovdi\@math.klte.hu
\bigskip
\hskip-\parindent {\sl V.P.~Rudko}
\newline
Department of Algebra,  University of Uzhgorod
\newline
88 000, Uzhgorod, Ukraine
\endaddress

\endtopmatter
\document
The first explicit description of the $\Bbb Z$-representations of
the noncyclic group $G$ of order $4$ was obtained by L.~Nazarova
\cite{4,5}. This description was significantly simplified and
clarified  after the discovery of  the connection between the
$\Bbb Z$-representations of the group $G$ and the representations
of a certain graph with five vertices  defined on  linear spaces
over the field of two elements \cite{2}.

Let $K=\Bbb Z$, $\Bbb Z_{(2)}$, or  $\Bbb Z_{2}$ be the ring of
rational integers, the localization of $\Bbb Z$ at the prime $2$
and the ring of $2$-adic integers, respectively. Let $G=\gp{a,
b}\cong C_2\times C_2$ be a noncyclic group of order $4$.

Let  $\frak{M}_n$ ($n>1$) be the set of all polynomials $f(x)$ of
degree $n$ over the field of two elements $\Bbb{F}_2$ which
satisfy one of the following conditions: \item{$\bullet$} $f(x)$
is irreducible over $\Bbb{F}_2$; \item{$\bullet$} $f(x)$ is a
power of a nonlinear irreducible  polynomial over $\Bbb{F}_2$.

Define an action of the group
$$
\split
Aut(G)=\gp{ \sigma_1,\sigma_2 \mid  &\sigma_1(a)=b, \quad \sigma_1(b)=a,\\
&\sigma_2(a)=a,\quad  \sigma_2(b)=ab,\quad  a,b\in G }\cong S_3
\endsplit
$$
on the set $\frak{M}_n$  by\quad  $
f^{\sigma_1}(x)=x^nf(x\m1)$, \quad  $f^{\sigma_2}(x)=f(x+1)$,
where $f(x)\in \frak{M}_n$ and $n=deg(f(x))$.

By $\frak{M}^\prime_n$ we denote the set of the orbits of the set
$\frak{M}_n$ under the action of the group $Aut(G)$. Let
$$
\matrix e_1=\frac{1}{4}(1+a)(1+b); &\quad
e_2=\frac{1}{4}(1-a)(1-b);\cr e_3=\frac{1}{4}(1+a)(1-b); &\quad
e_4=\frac{1}{4}(1-a)(1+b),\cr
\endmatrix\tag1
$$
be the four primitive idempotents of the  group algebra $FG$ over
a filed $F$  containing $K$. The corresponding irreducible
$K$-representations of these idempotents have  the following form:
$$
\matrix
\format \r\quad &\r\quad &\r\quad &\r\quad &\r\quad &\r\quad &\r\quad &\r\quad &\r\quad &\r\quad &\r\quad &\r\quad &\r\quad &\r\quad &\r\\
\chi_1:  & a & \mapsto & 1, & b &\mapsto &1; \quad& \chi_2:  & a
&\mapsto &-1, & b &\mapsto &-1;\cr \chi_3:  & a &\mapsto  & 1, & b
&\mapsto &-1;\quad& \chi_4:  & a &\mapsto  &-1, & b &\mapsto & 1.
\cr
\endmatrix\tag2
$$
The representations $\chi_i$  $(i=1, \ldots, 4)$ are realized in
the $KG$-module $KGe_i=Ke_i$ of the group algebra $FG$.

Consider  the following $K$-representations $\Gamma$ of the group
$G=\gp{a,b}=C_2\times C_2$ which were introduced in \cite{3}.  In
the following  we shell refer to them as ($\Bbb{L}$):
$$
\matrix
&&&&\\
\Delta_n:& a\mapsto &\left(\smallmatrix
E_n  & 0 &  0   &  E_n & 0    \\
     & 1 &  0   &  0   & 0    \\
     &   & -E_n &  0   & E_n  \\
     &   &      & -E_n & 0    \\
     &   &      &      & E_n  \\
\endsmallmatrix\right),
\quad \quad \quad \quad &b\mapsto &\left(\smallmatrix
1    &   0 &  0   &  0   &  0    \\
     & E_n &  0   &  0   &  E_n  \\
     &     & -E_n &  E_n &  0    \\
     &     &      &  E_n &  0    \\
     &     &      &      & -E_n  \\
\endsmallmatrix\right);
&&&&\\
&&&&\\
&&&&\\
&&&&\\
W_0:& a\mapsto & \left( \smallmatrix
1  &  1    &  0   &  1    \\
   & -1    &  0   &  0    \\
   &       &  1   &  0    \\
   &       &      & -1    \\
\endsmallmatrix\right),
& b\mapsto & \left( \smallmatrix
1  &  1    &  1   &  0    \\
   & -1    &  0   &  0    \\
   &       & -1   &  0    \\
   &       &      &  1    \\
\endsmallmatrix\right); \\
&&&&\\
&&&&\\
W_n:& a\mapsto& \left(\smallmatrix
D  & 0   & 0    &  0       &  0        \\
   & E_n & 0    &  0       &  \NL      \\
   &     & -E_n &  \NL     &  0        \\
   &     &      &  E_{n+1} &  0        \\
   &     &      &          & -E_{n+1}  \\
\endsmallmatrix\right),  &
b\mapsto& \left(\smallmatrix
D  &  0  &  0   &     S       & 0         \\
   & E_n &  0   &   \NP    & 0         \\
   &     & -E_n &     0       & \NP       \\
   &     &      &  -E_{n+1}   & 0         \\
   &     &      &             & E_{n+1}   \\
\endsmallmatrix\right);\\  &&&&\\
&&&&\\
&&&&\\
T_n:
    & a\mapsto& \left(\smallmatrix
E_{n+1}   &  0  &  E_{n+1}  &   0       \\
          &-E_n &  0        &   E_n     \\
          &     & -E_{n+1}  &   0       \\
          &     &           &   E_{n}   \\
                  \endsmallmatrix
                  \right),
    & b\mapsto& \left(\smallmatrix
1&         &  0  &  0        &   0     \\
&E_{n}     &  0  &  0        &   E_n       \\
&          &-E_n &  \NP      &   0     \\
&          &     & E_{n+1}   &   0       \\
&          &     &           &  -E_{n}   \\
                  \endsmallmatrix
                  \right);\\
&&&&\\
&&&&\\
\Delta_{n,1}:& a\mapsto& \left(\smallmatrix
E_{n}     &   0      & \NP       &  0        \\
          & -E_{n+1} &     0     &  E_{n+1}  \\
          &          &  -E_{n+1} &  0        \\
          &          &           &  E_{n+1}  \\
\endsmallmatrix\right),  &
b\mapsto& \left(\smallmatrix
E_{n}     &  0      &   0       & \NL   \\
          &-E_{n+1} &  E_{n+1}  &   0        \\
          &         & E_{n+1}   &   0        \\
          &         &           &   -E_{n+1} \\
\endsmallmatrix\right),\\
\endmatrix
$$
where $S=\left(\smallmatrix 0&\ldots&0& 1&1\\0&\ldots&0& 0&0\\
\endsmallmatrix\right)$,\quad  $D=\left(\smallmatrix 1&1\\0&-1\\
\endsmallmatrix\right)$.

The module $M$ is called  a module of the $K$-representation of
the group $G$ if $M$ is a $KG$-module with finite $KG$ basis.

\proclaim{Lemma 1} If $M$ is a module of the $K$-representation
$\Gamma$ of the group $G$ and $M$ does not contains regular direct
summands, then
$$
(1+a)(1+b)M\subset 2M. \tag3
$$
\endproclaim
\demo{Proof} Let $\overline{K}=K/2K$ be the field of $2$ elements.
From the description of the $\overline{K}$-representations of
group $G\cong C_2\times C_2$ (see \cite{1}) follows that if
$\overline{M}$ is a module of indecomposable
$\overline{K}$-representation of $G$ and
$\overline{M}\not=\overline{K}G$, then
$$
(1+a)(1+b)\overline{M}=0.\tag4
$$
Let ${M}$ be a module of an indecomposable ${K}$-representation of
$G$ and ${M}\not={K}G$. The description of the $K$-representations
of $G$ (see \cite{4,5}) shows that $\overline{M}=M/2M$ is a sum of
indecomposable $\overline{K}G$-modules, and it is different from
$\overline{K}G$. Hence from (4) follows (3), so the  proof is
complete.
\enddemo

Let $M$ be a module of the $K$-representation of the group $G\cong
C_2\times C_2$, which satisfies the condition of Lemma 1. Put
$M'=e_1M+e_2M+e_3M+e_4M\subset FM$, where $V_0=M'/M$, $M_i=e_iM+M$
and  $V_i=M_i/M$ ($i=1,2,3,4$). Clearly, by Lemma 1 we get that
$2M'\subset M$, $2M_i\subset M$ and $V_0, V_i$ ($i=1,2,3,4$) are
vector spaces over the  field $\Bbb{F}_2=K/2K$. These spaces form
a diagram
$$
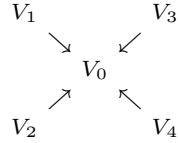
\smallmatrix
V_1& & & &V_3\\
&\searrow & &\swarrow &  \\
& &V_0 & &  \\
&\nearrow & &\nwarrow &  \\
V_2& & & &V_4\\
\endsmallmatrix
$$
in which the arrows correspond to homomorphisms generated by the
embeddings of $M_i$ to $M'$ $(i=1, \ldots, 4)$. Thus, the 5-tuple
of linear spaces $(V_0; V_1, V_2, V_3, V_4)$ is a representation
of the following graph on five points:
$$
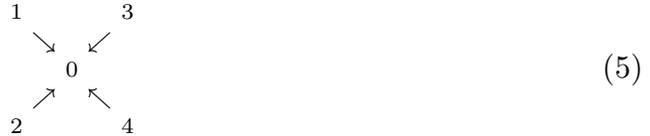
\smallmatrix
1& & & &3\\
&\searrow & &\swarrow &  \\
& &0 & &  \\
&\nearrow & &\nwarrow &  \\
2& & & &4\\
\endsmallmatrix
\tag5
$$
The following quadratic form
$$
B(x)=x_0^2+x_1^2+x_2^2+x_3^2+x_4^2-x_0(x_1+x_2+x_3+x_4)
$$
is defined on $5$-dimension vectors corresponding  to the graph
(5), where $x=(x_0; x_1, x_2, x_3, x_4) \in {\Bbb Z}^5$.

A vector $x=(x_0; x_1,x_2,x_3,x_4)\in {\Bbb Z}^5$ is called a {\it
positive root} of the form $B$ if $B(x)=1$ and $x_i\geq 0$, where
$i=0,\ldots, 4$. The {\it dimension vector } $d(\Gamma)=(d_0; d_1,
d_2, d_3, d_4)$ of a $K$-representation $\Gamma$ of the group $G$
is denoted by $d_i=dim_{\Bbb{F}_2}V_i$, where $V(\Gamma)=(V_0;
V_1, V_2, V_3, V_4)$ is the 5-tuple of linear spaces corresponding
to the representation $\Gamma$.

We know (see \cite{2}),  that if $\Gamma$ is an indecomposable
$K$-representation of the group $G$, which is different from the
regular and irreducible representations, then either
$B(d(\Gamma))=0$ or $d(\Gamma)$ is a positive root of the form
$B(x)$ (i.e. $B(d(\Gamma))=1$). In addition, $d_0>0$ and
$d_1+d_2+d_3+d_4>1$ for any positive root $d=(d_0; d_1, d_2, d_3,
d_4)$. Moreover, for every positive root $d=(d_0; d_1, d_2, d_3,
d_4)$ such that $d_0>1$ and $d_1+d_2+d_3+d_4>1$ there exists only
one indecomposable $K$-representation of the group $G$, such that
the degree $\Gamma$ is equal to $d_1+d_2+d_3+d_4$.

Clearly, the one parameter family of indecomposable
$K$-representations of the group $G$ of degree $4m$ corresponds to
the vector $d=(2n; n, n, n, n)$ (of course, $B(d)=0$). The
positive roots of the form $B$ are listed in the following table:
\bigskip
{ \eightpoint{
\centerline{\vbox{\halign{\strut\offinterlineskip\vrule\quad\hfill
$#$\hfill&\quad
\vrule\quad\hfill$#$\hfill&\quad\vrule\quad\hfill$#$\hfill\quad\vrule\cr
\noalign{\hrule} &&\cr \noalign{\vskip -5pt}
m=d_1+d_2+d_3+d_4&d=(d_0; d_1, d_2, d_3, d_4) &t(m)\cr
\noalign{\vskip -3pt} &&\cr \noalign{\hrule} &&\cr \noalign{\vskip
-3pt}
4n& (2n+1; n, n, n, n)&  1\cr
\noalign{\vskip -3pt}
n\geq 1&&\cr
\noalign{\hrule}
4n& (2n-1; n, n, n, n)&  1\cr
\noalign{\vskip -3pt}
n\geq 1&&\cr
\noalign{\hrule}
4n+2& (2n+1; n+1, n+1, n, n)&  6\cr
\noalign{\vskip -3pt}
n\geq 1&&\cr
\noalign{\hrule}
4n+1& (2n; n+1, n, n, n)&  4\cr
\noalign{\vskip -3pt}
n\geq 1&&\cr
\noalign{\hrule}
4n+2& (2n+1; n+1, n, n, n)&  4\cr
\noalign{\vskip -3pt}
n\geq 1&&\cr
\noalign{\hrule}
4n+3& (2n+2; n+1, n+1, n+1, n)&  4\cr
\noalign{\vskip -3pt}
n\geq 1&&\cr
\noalign{\hrule}
4n+3& (2n+1; n+1, n+1, n+1, n)&  4\cr \noalign{\vskip -3pt} n\geq
1&&\cr \noalign{\hrule} }}}}} \centerline{Table 1.}
\bigskip
\noindent where $t(m)$ denotes the number of roots with prescribed
$m$ and $d_0$. These roots are obtained from the  root $d$ by
permutations of its components $d_1, d_2, d_3, d_4$.

First we describe the indecomposable $K$-representations $\Gamma$
of $G$ with  common dimension vector $d(\Gamma)=(2n; n, n, n, n)$
(we recall that $B(d(\Gamma))=0$).

We will consider each  polynomial $f(x) \in \Bbb{F}_2[x]$ as a
polynomial over $K$, replacing $0,1 \in \Bbb{F}_2$ by  $0,1 \in
K$, respectively. Define  $ \widetilde{ f(x)}=\left(\smallmatrix
0&      &       &\alpha_0    \\
1&      &      &\alpha_1    \\
 &\ddots&       & \vdots      \\
0&      &   1   &\alpha_{n-1}\\
\endsmallmatrix\right)$ with the help of the corresponding
matrix of
$f(x)=x^n-\alpha_{n-1}x^{n-1}-\cdots-\alpha_1x-\alpha_0\in
\Bbb{F}_2[x]$.

Each polynomial $f(x) \in \frak{M}_n \cup \{ x^n, (x+1)^n \}$,
where $ \frak{M}_1 = \varnothing $, defines the following
indecomposable $K$-representation of the group $G$:

$$
\Delta_f: a\mapsto \left(\smallmatrix
E_n  &   0  & U_{11}  &  0     \\
     & -E_n &  0      & U_{22} \\
     &      &  -E_n   &  0     \\
     &      &         &  E_n   \\
\endsmallmatrix\right);
\quad \quad \quad b\mapsto \left(\smallmatrix
 E_n &  0   &  0    &  U_{12}  \\
     & -E_n & U_{21}&  0    \\
     &      & E_n  &  0    \\
     &      &       & -E_n  \\
\endsmallmatrix\right),
$$
where only one of the matrices $U_{ij}$ ($1\leq i,j\leq 2$) is
equal to $ \widetilde{f(x)}$, and the three other $U_{ij}$ are
equal to $E_n$.

Suppose  that $\frak F(f(x))$ is exactly  that representation
$\Delta_f$ for which  $U_{12}= \widetilde{f(x)}$. For every
$K$-representation $\Gamma$ of the group $G$ and for an
automorphism $\phi\in Aut(G)$ we define a $K$-representation
$\Gamma^\phi$ of the group $G$ by
$
\Gamma^\phi(g)=\Gamma(\phi(g)), $ where $g\in G$. The
$K$-representations   $\Gamma$ and $\Gamma^\phi$ are called {\it
conjugate}.

\proclaim{Lemma  2} Let $\Gamma$ be either an indecomposable
$K$-representation from the list $(\Bbb{L})$ or its tensor product
with the character $\chi_2$ from (2). Then we have:
\bigskip

{ \eightpoint{ \centerline{\vbox{\halign{\strut\offinterlineskip
\vrule\quad\hfill $#$\hfill&\quad \vrule\quad\hfill$#$\hfill\quad
\vrule\cr \noalign{\hrule} &\cr \noalign{}
\Gamma & d(\Gamma)=(d_0; d_1, d_2, d_3, d_4) \cr \noalign{\hrule}
\Delta_{n,1} & (2n+1; n,   n+1, n+1, n+1) \cr \noalign{\hrule}
\Delta_{n}^*& (2n+1; n+1, n, n, n)\cr \noalign{\hrule}
\Delta_{n}& (2n+1; n+1, n, n, n)\cr \noalign{\hrule}
W_n & (2n+1; n+1, n+1, n+1, n+1)\cr \noalign{\hrule}
T_n & (2n+1; n+1, n,   n, n+1)\cr \noalign{\hrule}
\Delta_{n,1}^*& (2n+2; n,n+1,   n+1, n+1)\cr \noalign{\hrule}
W_n^*& (2n+3; n+1,n+1,   n+1, n+1)\cr \noalign{\hrule}
\Delta_{n}\otimes\chi_2 &(2n;n,n+1,n,n)\cr \noalign{\hrule}
(\Delta_{n}\otimes\chi_2)^*& (2n;n,n+1,n,n)\cr \noalign{\hrule}
\Delta_{n,1}\otimes\chi_2 & (2n+1; n+1,n,   n+1, n+1)\cr
\noalign{\hrule}
(\Delta_{n,1}\otimes\chi_2)^*& (2n+2; n+1,n, n+1, n+1)\cr
\noalign{\hrule}
T_n\otimes\chi_2&(2n+1; n, n+1,   n, n+1) \cr \noalign{\hrule}
 }}}}}
 \centerline{Table  2.}
\endproclaim

\demo{Proof} We prove it for $\Gamma=\Delta_{n,1}$. Let
$$
\{v_1,\ldots,v_n,u_1,\ldots,u_{n+1},w_1,\ldots,w_{n+1},f_1,\ldots,f_{n+1}\}
$$
be a basis of the module $M$ of the $K$-representation $\Delta_{n,1}$. Clearly,
$$
\split
&e_1v_i=\textstyle{\frac{1}{4}}(1+a)(1+b)v_i=v_i;\quad\quad\quad\quad
e_1u_i=0;\\
&e_1w_i=\textstyle{\frac{1}{4}}(1+a)(2w_i+u_i)=\cases \textstyle{\frac{1}{2}}v_i, \quad &
i\not=n+1;\\ 0, \quad &
i=n+1;\\
\endcases\\
& e_1f_i=\textstyle{\frac{1}{4}}(1+b)(2f_i+u_i)=\cases \textstyle{\frac{1}{2}}v_{i-1},
\quad & i=2,\ldots,n+1;\\ 0, \quad &
i=1;\\
\endcases\\
& e_2v_i=\textstyle{\frac{1}{4}}(1-a)(1-b)v_i=0;\quad\quad
  e_2u_i=u_i;\\
& e_2w_i=\textstyle{\frac{1}{4}}(1-a)u_i=
\textstyle{\frac{1}{2}}u_i;\quad\quad\quad\quad
  e_2f_i=\textstyle{\frac{1}{4}}(1-b)(-u_i)=- \textstyle{\frac{1}{2}}u_i;\\
\endsplit
$$$$
\split & e_3v_i=\textstyle{\frac{1}{4}}(1+a)(1-b)v_i=0;\quad\quad
  e_3u_i=0;\quad
 e_3w_i=-\textstyle{\frac{1}{4}}(1+a)u_i=0; \\
&  e_3f_i=\textstyle{\frac{1}{4}}(1-b)(2f_i+u_i)=
\cases
f_i+\textstyle{\frac{1}{2}}(u_i-v_{i-1}),\quad & i=2,\ldots,n+1;\\
f_1+\textstyle{\frac{1}{2}}u_1,\quad & i=1;\\
\endcases
\endsplit
$$$$
\split & e_4v_i=\textstyle{\frac{1}{4}}(1-a)(1+b)v_i=0;\quad\quad
  e_4u_i=0;\\
&
e_4w_i=-\textstyle{\frac{1}{4}}(1-b)u_i=-\textstyle{\textstyle{\frac{1}{2}}}u_i;
\quad
  e_4f_i=-\textstyle{\frac{1}{4}}(1+b)u_i=0.\\
\endsplit
$$
Since $V_i=(e_iM+M)/M$ we have that
$$
\split V_1=&\gp{\textstyle{\frac{1}{2}}v_{i}+M\mid
i=1,\ldots,n};\quad
V_2=\gp{\textstyle{\frac{1}{2}}u_{i}+M\mid i=1,\ldots,n+1};\\
V_3=&\gp{\textstyle{\frac{1}{2}}u_{1}+M,\quad \textstyle{\frac{1}{2}}(u_i-v_{i-1})+M\mid
i=2,\ldots,n+1};\\
V_4=&\gp{-\textstyle{\frac{1}{2}}u_i+M\mid i=1,\ldots,n+1}.\\
\endsplit
$$
Finally, it is easy to check that
$$
V_0=\gp{\textstyle{\frac{1}{2}}u_i+M,\;
\textstyle{\frac{1}{2}}v_j+M\mid i=1,\ldots,n+1; \; j=1,\ldots,n}
$$
and  $ d(\Delta_{n,1}) =(2n+1; n, n+1, n+1, n+1)$, so the proof is
complete. The proofs for other representations  $\Gamma$ are
similar to this one. The lemma is proved.
\enddemo

\proclaim{Lemma 3} Let $K\in\{ \Bbb Z, \Bbb Z_{(2)},\Bbb Z_{2}\}$.
Up to conjugacy, all indecomposable $K$-repre\-sen\-tations
$\Gamma$ of $G\cong C_2\times C_2$, with  common dimension vector
$d(\Gamma)=(2n; n, n, n, n)$ belong to the either $\frak{F}(f(x))$
\quad ($f(x)\in \frak{M}_n\cup\{x^n, (x+1)^n\}$) or to
$\frak{F}'(f(x))=\frak{F}(x^n)\otimes \chi_2$, where $\chi_2$ is
from (2).

Up to equivalence, all indecomposable $K$-representations $\Gamma$
of $G$ with  $d(\Gamma)=(2n; n, n, n, n)$ belong to  one of the
following types: \item{$\bullet$} $\frak{F}(f(x))$, \qquad
$f(x)\in \frak{M}_n\cup\{x^n, (x+1)^n\}$; \item{$\bullet$}
$(\frak{F}(x^n))^{\sigma_1}$; \quad \item{$\bullet$}
$\frak{F}'(x^n)$;\quad \item{$\bullet$}
$(\frak{F}'(x^n))^{\sigma_1}$;\quad \item{$\bullet$}
$(\frak{F}'(x^n))^{\sigma_2}$, \quad where $\sigma_1,\sigma_2\in
Aut(G)$.
\endproclaim

\demo{Proof} The analysis of well known Nazarova's results on the
representations of the group  $G$ (see \cite{3,4}) shows that an
indecomposable $K$-representation $\Gamma$ of $G$, such that
$d(\Gamma)=(2n; n, n, n, n)$, is equivalent to one of the
representations $\Delta_f$,\quad    where  $f(x)\in
\frak{M}_n\cup\{x^n, (x+1)^n\}$,\quad   except the
$K$-representation $(\frak{F}'(x))^{\sigma_2}$, which is
equivalent to the following representation:
$$
a\mapsto \left(\smallmatrix
L_1&  0     & 0   &  0     &  0        & 0 \\
   & E_{n-1}& 0   &  0     &E_{n-1}    & 0 \\
   &        &E_{n-1}& E_{n-1}&  0     & 0 \\
   &        &        &-E_{n-1}&  0     & 0 \\
   &        &        &        & E_{n-1}&0\\
   &        &        &        &        &-L_1\\
\endsmallmatrix\right);
\quad b\mapsto \left(\smallmatrix
L_1&  0     & 0   &        &  L_2      & 0 \\
   & E_{n-1}& 0   &  0     &J_{n-1}(1) & L_3 \\
   &        &-E_{n-1}&-E_{n-1}&  0     & 0 \\
   &        &        & E_{n-1}&  0     & 0 \\
   &        &        &        & E_{n-1}&0\\
   &        &        &        &        &L_1\\
\endsmallmatrix\right),
$$
where $J_{k}(1)$ is a Jordan block of dimension $k$ with ones  in
the main diagonal,
\newline
$L_1=\left(\smallmatrix
1&1\\
0&-1\\
\endsmallmatrix\right)$,
\quad
$L_2=\left(\smallmatrix
1&\ldots &0\\
0&\ldots &0\\
\endsmallmatrix\right)$,
\quad
$L_3=\left(\smallmatrix
0      & 0     \\
\vdots & \vdots\\
0      & 1     \\
\endsmallmatrix\right)$. So the lemma is proved.
\enddemo

\newpage

In the next Theorem among other $K$-representations of the group
$G$ we list  those indecomposable $K$-representations $\Gamma$
such that $d(\Gamma)$ is a root of the form $B(x)$ (i.e.
$B(d(\Gamma))=1$). Analyzing Nazarova's classification \cite{4,5},
we can see that any such obtained $K$-representation $\Gamma$ can
be derived from $\Omega$ and $\chi_i$, (where $\Omega$ is from
($\Bbb{L}$), $\chi_i$ is from (2)) by one of the following
operations:

\itemitem{$\bullet$} taking its conjugate $\Omega^\phi$, where \;
$\phi\in Aut(G)$; \itemitem{$\bullet$} taking the contragradient
$K$-representation $\Omega^*$ ($\Omega^*(g)=\Omega^T(g^{-1})$, \;
$g\in G$); \itemitem{$\bullet$} taking the tensor product
$\Omega\otimes\chi_i$. \smallskip

Our main result reads as follows:
 \proclaim {Theorem} Let $K=\Bbb Z$, $\Bbb Z_{(2)}$, or $\Bbb
Z_{2}$ be the ring of rational integers, the localization of $\Bbb
Z$ at the prime $2$ and the ring of $2$-adic integers,
respectively. All indecomposable subgroups in the group $GL(m,K)$
which are isomorphic to the group $G=\gp{a,b}\cong C_2\times C_2$
can be described, up to conjugacy, by groups which are generated
by matrices $\Gamma(a)$ and $\Gamma(b)$, where $\Gamma$ is a
$K$-representation of $G$ from the following table:
\endproclaim

{ \eightpoint{
\centerline{\vbox{\halign{\strut\offinterlineskip\vrule\quad\hfill
$#$\hfill&\quad
\vrule\quad\hfill$#$\hfill&\quad\vrule\quad\hfill$#$\hfill&\quad
\vrule\quad\hfill$#$\hfill&\quad\vrule\quad\hfill$#$\hfill\quad\vrule\cr
\noalign{\hrule} &&&&\cr \noalign{\vskip -5pt}
m&\Gamma&d(\Gamma)&St(\Gamma)&S(\Gamma)\cr \noalign{\vskip -3pt}
&&&&\cr \noalign{\hrule} &&&&\cr \noalign{\vskip -3pt}
4n+1& \Delta_n& (2n+1; n+1, n, n, n)& Aut(G) & 1\cr
\noalign{\vskip -3pt} n\geq 1&&&&\cr \noalign{\hrule}
4n+1& \Delta_n^*& (2n+1; n+1, n, n, n)& Aut(G) & 1\cr
\noalign{\vskip -3pt} n\geq 1&&&&\cr \noalign{\hrule}
4n+1& \Delta_n\otimes \chi_2& (2n+1; n, n+1, n, n)& \gp{\sigma_1}
& 3\cr \noalign{\vskip -3pt} n\geq 1&&&&\cr \noalign{\hrule}
4n+1& \Delta_n^*\otimes \chi_2& (2n; n, n+1, n, n)& \gp{\sigma_1}
& 3\cr \noalign{\vskip -3pt} n\geq 1&&&&\cr \noalign{\hrule}
4n+3& \Delta_{n,1}& (2n+1; n, n+1, n+1, n+1)& Aut(G) & 1\cr
\noalign{\vskip -3pt} n\geq 0&&&&\cr \noalign{\hrule}
4n+3& \Delta_{n,1}^*& (2n+2; n, n+1, n+1, n+1)& Aut(G) & 1\cr
\noalign{\vskip -3pt} n\geq 0&&&&\cr \noalign{\hrule}
4n+3& \Delta_{n,1}\otimes\chi_2& (2n+1; n+1, n, n+1, n+1)&
\gp{\sigma_1} & 3\cr \noalign{\vskip -3pt} n\geq 0&&&&\cr
\noalign{\hrule}
4n+3& \Delta_{n,1}^*\otimes\chi_2& (2n+2; n+1, n, n+1, n+1)&
\gp{\sigma_1} & 3\cr \noalign{\vskip -3pt} n\geq 0&&&&\cr
\noalign{\hrule}
4n+2& T_n& (2n+1; n+1, n, n+1,n)& \gp{\sigma_1\sigma_2\sigma_1} &
3\cr \noalign{\vskip -3pt} n\geq 1&&&&\cr \noalign{\hrule}
4n+2& T_n\otimes\chi_2& (2n+1; n, n+1, n, n+1)&
\gp{\sigma_1\sigma_2\sigma_1} & 3\cr \noalign{\vskip -3pt} n\geq
1&&&&\cr \noalign{\hrule}
4n+4& W_n& (2n+1; n+1, n+1, n+1, n+1)& Aut(G) & 1\cr
\noalign{\vskip -3pt} n\geq 0&&&&\cr \noalign{\hrule}
4n+4& W_n^*& (2n+3; n+1, n+1, n+1, n+1)& Aut(G) & 1\cr
\noalign{\vskip -3pt} n\geq 0&&&&\cr \noalign{\hrule}
4n& \frak{F}(f(x)),& & St[f(x)] & \frac{6}{|St[f(x)]|}\cr
\noalign{\vskip -3pt} n\geq 1&\quad (f(x)\in \frak{M}_n')&&&\cr
\noalign{\hrule}
4n& \frak{F}(x^n)& &  \gp{\sigma_2\sigma_1\sigma_2} & 3\cr
\noalign{\vskip -3pt} n\geq 1&&&&\cr \noalign{\hrule}
4n& \frak{F}'(x^n)=& &  \gp{\sigma_2\sigma_1\sigma_2} & 3\cr
\noalign{\vskip -3pt} n\geq 1&\frak{F}(x^n)\otimes\chi_2&&&\cr
\noalign{\hrule}
4&\Gamma \text{ is the regular rep.}& &   & 1\cr \noalign{\hrule}
}}}}} \centerline{Table 3}


\noindent {\it Here  $St(\Gamma)=\{ \phi\in Aut(G)\; |\;
\Gamma^\phi\;  \text{ is equivalent to }\; \Gamma \}$ is the
stabilizer of $\Gamma$ in $Aut(G)$,\quad
$S(\Gamma)=\frac{|Aut(G)|}{|St(\Gamma)|}$ is the number of
non-equivalent $K$-representations which are conjugate to
$\Gamma$, \quad  $\Gamma^*$ is the contragradient
$K$-representation of  $\Gamma$, defined by
$\Gamma^*(g)=\Gamma^T(g\m1)$\quad $(g\in G)$}.

\subhead Proof of the Theorem
\endsubhead
Let $G_1$ and $G_2$ be  subgroups in $GL(m,K)$ isomorphic to
$G=\gp{a,b}\cong C_2\times C_2$ with isomorphisms $\Gamma_i:
G\mapsto G_i$, ($i=1,2$). Then $\Gamma_i$ is a faithful
$K$-representation of $G$. By definition, the $K$-representation
$\Gamma_1$ is conjugate to the $K$-representation $\Gamma_2$, if
there are exist $C\in GL(m,K)$, $\phi\in Aut(G)$, such that
$C\m1\Gamma_1(\phi(g))C=\Gamma_2(g)$, where $g\in G$. In other
words, $K$-representations $\Gamma_1$ and $\Gamma_2$ of the group
$G$ are conjugate, if the groups $\Gamma_1(G)$ and $\Gamma_2(G)$
are conjugate in $GL(m,K)$. This implies  the "if" part of the
Theorem.

Now suppose $G_1$ and $G_2$ are conjugate in $GL(m,K)$ and $\tau:
G_1 \mapsto G_2$ is an isomorphism such that  $\tau(g)=C\m1gC$ for
some matrix $C\in GL(m,K)$, where  $g\in G$. Obviously,  the map
$$
\phi:G @>{\Gamma_1}>>G_1@>{\tau}>>  G_2 @>{\Gamma_2\m1}>>G
$$
is an automorphism of $G$ and
$$
\Gamma_2^\phi(g)=\Gamma_2(\Gamma_2\m1\tau\Gamma_1(g))=\tau(\Gamma_1(g))=
C\m1\Gamma_1(g)C,
$$
where $g\in G$, i.e. the $K$-representations $\Gamma_1$ and
$\Gamma_2$ are conjugate.

Let $\Gamma$ be a faithful $K$-representation of $G$, let
$St(\Gamma)$ be the stabilizer of $\Gamma$ in $Aut(G)$ and let
$\{g_1,\ldots,g_t\}$ be a  system of representatives of the cosets
of $Aut(G)$ by $St(\Gamma)$. Then $\{\Gamma^{g_i}\; \mid \;
i=1,\ldots,t\}$ is the set of the pairwise nonequivalent
$K$-representations of $G$, which are conjugated with $\Gamma$.

If $\Gamma$ is an indecomposable $K$-representation of $G$ and
$d(\Gamma)$ is a root of the form $B(x)$, then $d(\Gamma^\phi)$
($\phi\in Aut(G)$) is also a root of the form $B(x)$. The vector
$d(\Gamma^\phi)$ can be obtained from the vector $d(\Gamma)=(d_0;
d_1, d_2, d_3, d_4)$ by a permutation of the last four components
of $d(\Gamma)$. It follows, that $S(\Gamma)$ (see Table 1) is the
number of roots, which we obtained  from $d(\Gamma)$ by
permutation of some components in $d(\Gamma)$.

To finish the proof, we must  compare the number of
$K$-representations of degree $m$ in Table 1 to the sum of all
$t(m)$ for a given $m$ in the Table 3, and then we can apply Lemma
3. So the lemma is proved.
\bigskip

\proclaim {Remark} Let us substitute  each $K$-representation
$\Gamma$ of $G$ in  Table 2 by  $\Gamma^\sigma$. Here $\sigma$
denotes  a coset  of $Aut(G)$ by $St(\Gamma)$. In this manner  we
obtain all faithful indecomposable $K$-representations of $G$ up
to equivalence.
\endproclaim
The proof of this {\it Remark} follows from the analysis of the
classification of the representations described by L.~Nazarova
\cite{3,4}.


\enddocument